\documentclass[preprint,11pt]{article}
\usepackage{a4,amsmath,amsthm,amsfonts,amssymb,hyperref}
\usepackage{graphicx,psfrag}

\title{On the residual finiteness of outer automorphisms of relatively hyperbolic groups}
\author{V. Metaftsis and M. Sykiotis}

\makeatletter

\newtheorem{thm}{Theorem}[section]
\newtheorem{lem}[thm]{Lemma}

\newtheorem{claim}{Claim}
\theoremstyle{definition}

\theoremstyle{remark}
\newtheorem{rem}[thm]{Remark}

\def\classification{\@ifnextchar [{\@xfootnotetext}%
{\begingroup\let\protect\noexpand\xdef\@thefnmark{}
\endgroup\@footnotetext}}
\def\keywords{\@ifnextchar [{\@xfootnotetext}%
{\begingroup\let\protect\noexpand\xdef\@thefnmark{}
\endgroup\@footnotetext}}

\begin{document}
\classification {2000 {\sl Mathematics Subject Classification.}
20F65, 20F67, 20E26, 20E36, 20F28.} \keywords {{\sl Keywords and
phrases.} Relatively hyperbolic groups, Residually finite, Graphs
of groups, Limit groups, Outer automorphisms.} \maketitle

\begin{abstract} We show that every virtually torsion-free subgroup of the outer automorphism group
of a conjugacy separable relatively hyperbolic group is residually
finite. As a direct consequence, we obtain that the outer
automorphism group of a limit group is residually finite.
\end{abstract}

\section{Introduction}
Relatively hyperbolic groups were introduced by Gromov in
\cite{Gro}, in order to generalize notions such as the fundamental
group of a complete, non-compact, finite volume hyperbolic
manifold and to give a hyperbolic version of small cancellation
theory over free groups by adopting the geometric language of
manifolds with cusps. This notion has been developed by several
authors and, in particular, various characterizations of
relatively hyperbolic groups have been given (see \cite{Bo,Os} and
\cite{DS} and references therein). We should mention here that
Farb \cite{Fa} introduced a weaker notion of relative
hyperbolicity for groups, using constructions on their Cayley
graphs, as well as the Bounded Coset Penetration property, an
additional condition which makes his definition equivalent to the
other definitions.
\par
We recall here one of Bowditch's equivalent definitions (in the
case of infinite ``peripheral" subgroups). A finitely generated
group $G$ is \emph{hyperbolic relative to a family of finitely
generated subgroups} $\mathcal{G}$ if $G$ admits a proper,
discontinuous and isometric action on a proper, hyperbolic path
metric space $X$ such that $G$ acts on the ideal boundary of $X$
as a geometrically finite convergence group and the elements of
$\mathcal G$ are the maximal parabolic subgroups of $G$.

Besides the fundamental groups of hyperbolic manifolds of finite
volume, examples of relatively hyperbolic groups are fundamental
groups of finite graphs of finitely generated groups with finite
edge groups, which are hyperbolic relative to the family of
infinite vertex groups (which may be empty, in which case the
group is hyperbolic), since their action on the Bass-Serre tree
satisfies Definition 2 in \cite{Bo}.

Another example of relatively hyperbolic groups are limit groups.
The notion of a limit group was introduced by Sela \cite{S2,S} in
his solution to Tarski's problem for free groups. As it turned
out, the family of limit groups coincides with that of finitely
generated fully residually free groups first introduced by
Baumslag \cite{Ba}, and extensively studied by Kharlampovich and
Myasnikov \cite{KM1,KM2}. In \cite{Da}, Dahmani showed that limit
groups are hyperbolic relative to their maximal non-cyclic abelian
subgroups (see also \cite{Al}). Note that each group is relatively
hyperbolic to itself. So, from now on, in order to avoid this
trivial situation, \emph{we assume that all relatively hyperbolic
groups properly contain the corresponding maximal parabolic
subgroups}.
\par
In \cite{MS}, it was proved that the outer automorphism group of a
conjugacy separable hyperbolic group is residually finite. This is
a far-reaching generalization of a classical result of Grossman
\cite{Gr}, which states that the mapping class group of a closed
orientable surface is residually finite. The purpose of this note
is to prove the following generalization for relatively hyperbolic
groups.

\begin{thm}\label{thm:1}
Let $G$ be a conjugacy separable, relatively hyperbolic group.
Then each virtually torsion-free subgroup of the outer
automorphism group ${\rm Out}(G)$ of $G$ is residually finite.
\end{thm}
As an application, we obtain:

\begin{thm}\label{thm:2}
Let $G$ be the fundamental group of a finite graph of groups, such
that each edge group is finite and each vertex group is
polycyclic-by-finite. Then ${\rm Out} (G)$ is residually finite.
\end{thm}

Guirardel and Levitt \cite{GL} showed that the outer automorphism
group of a limit group is virtually torsion-free. More recently,
Chagas and Zalesskii \cite{CZ} have shown that limit groups are
conjugacy separable. Therefore, from Theorem \ref{thm:1}, we
immediately deduce the following result.

\begin{thm} The outer automorphism group of a limit group is residually finite.
\end{thm}

\section{Proofs of the main results}

A group $G$ is \emph{conjugacy separable} if for any two
non-conjugate elements $x$ and $y$ of $G$, there is a finite
homomorphic image of $G$ in which the images of $x$ and $y$ are
not conjugate. An automorphism $f$ of a group $G$ is called
\emph{conjugating} if $f(g)$ is conjugate to $g$ for each $g\in
G$. The conjugating automorphisms of a group $G$ form a subgroup
of ${\rm Aut}(G)$, which we denote by ${\rm Conj}(G)$. Clearly,
${\rm Conj}(G)$ is a normal subgroup of ${\rm Aut}(G)$ containing
the inner automorphism group ${\rm Inn}(G)$ of $G$. The importance
of this notion to the study of residual properties of the outer
automorphism group of $G$, arises from two facts. The first is
that if $G$ is finitely generated and conjugacy separable, then
the quotient group ${\rm Aut}(G)/{\rm Conj}(G)$ is residually
finite (see \cite[Lemma 2.1]{MS}). The second is the following
short exact sequence
\begin{equation}\label{eq:1} 1\rightarrow {\rm Conj}(G)/{\rm Inn}(G)\hookrightarrow
{\rm Aut}(G)/{\rm Inn}(G)\rightarrow {\rm Aut}(G)/{\rm
Conj}(G)\rightarrow 1.\end{equation}

Thus to prove Theorem \ref{thm:1}, it suffices to show that if $G$
is relatively hyperbolic, then the quotient ${\rm Conj}(G)/{\rm
Inn}(G)$ is finite.

We will need the following lemma whose a more general version in
the case of projections on quasiconvex subspaces can be found in
\cite[Proposition 2.1, Chapter 10]{CDP}.

\begin{lem} \label{lem:inequalities} Let $(X,d)$ be a $\delta$-hyperbolic metric space, let $g$
be an isometry of $X$ and let $N$ be a positive real number such
that the set $Y=\{y\in X:\,d(y,gy)\leq N\}$ is nonempty. Given a
point $x\in X$ and a positive number $M$, choose $y\in Y$ with
$d(x,y)\leq d(x,Y)+M$. Then either $d(x,gx)\geq
2d(x,y)+d(y,gy)-2(3\delta+2M)$ or $d(y,gy)\leq 3\delta+2M$.
\end{lem}

\begin{proof} We consider a geodesic trianlge with vertices $x$,
$y$ and $gy$ (see Figure 1). Let $z$ and $w$ be the points on the
geodesics $[x,y]$ and $[y,gy]$, respectively, which are at
distance $\alpha$ from $y$, where $\alpha$ is the Gromov product
of $x$ and $gy$ with respect to $y$. We first note that $w\in Y$.
Indeed, $d(w,gw)\leq d(w,gy)+d(gy,gw)=d(w,gy)+d(y,w)=d(y,gy)\leq
N$. Since $w\in Y$, we have $d(w,x)\geq d(x,Y)\geq d(x,y)-M$ and
hence $d(x,z)+d(z,y)=d(x,y)\leq d(w,x)+M\leq d(x,z)+\delta + M$.
It follows that $\alpha =d(z,y)\leq \delta + M $. Therefore
$d(x,gy)=d(x,y)+d(y,gy)-2\alpha$, where $\alpha\leq \delta + M$.
Similarly one can show that $d(y,gx)=d(x,y)+d(y,gy)-2\beta$, where
$\beta\leq \delta + M$. Now we turn our attention to a geodesic
quadrilateral with vertices $x$, $y$, $gy$ and $gx$. There are two
cases to consider, as shown in the following figure.

\hspace*{-1.5cm} \psfrag{xi}{$x$} \psfrag{yi}{$y$}
\psfrag{fi(g)xi}{$gx$} \psfrag{fi(g)yi}{$gy$} \psfrag{z}{$z$}
\psfrag{a}{$\alpha$} \psfrag{w}{$w$} \psfrag{Case 1.}{Case 1}
\psfrag{Case 2.}{Case 2}
\includegraphics[scale=0.45]{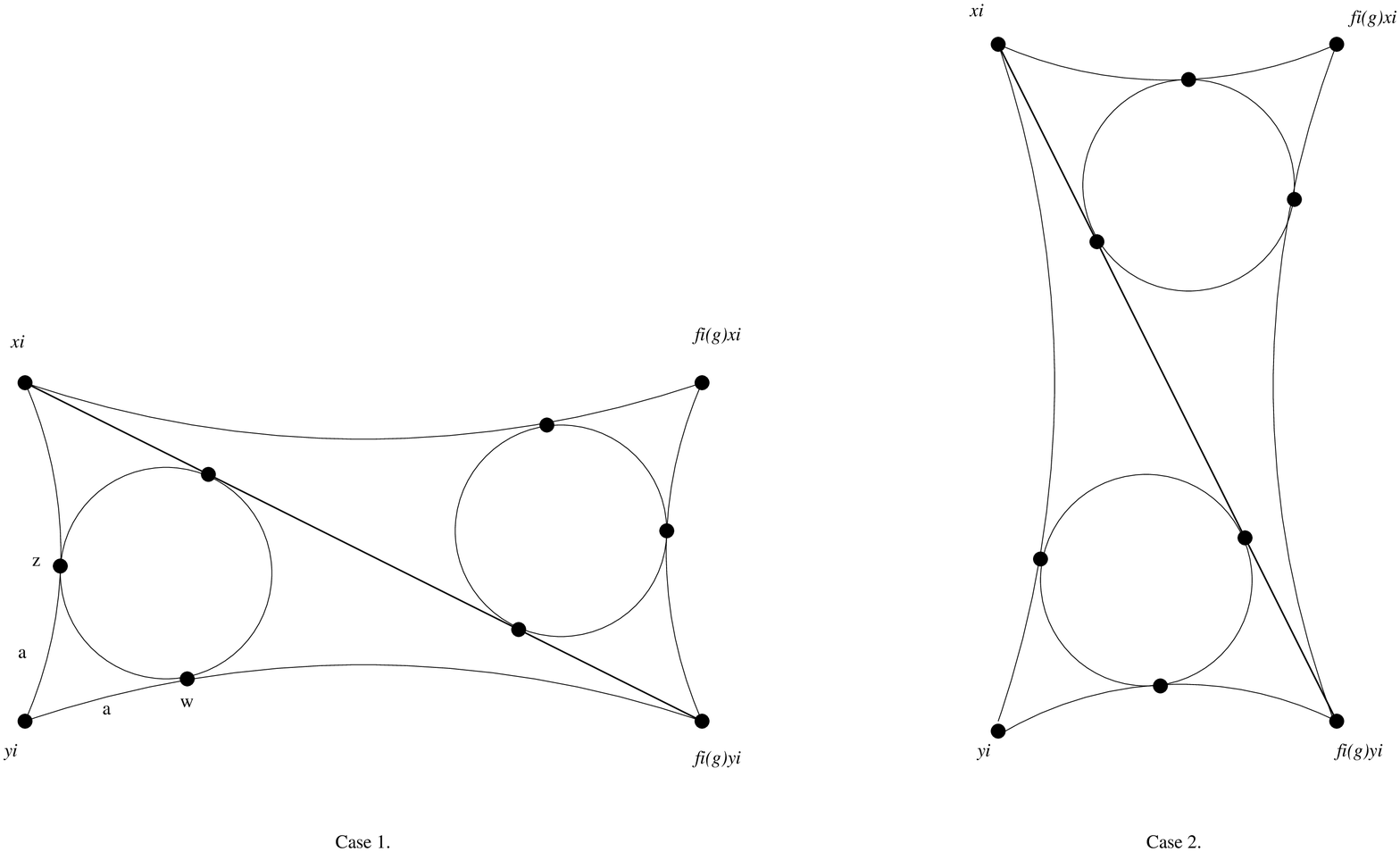}
\begin{figure}[h]  
\caption{The two cases of Lemma \ref{lem:inequalities}}
\end{figure}
\vspace*{.5cm}

In the first case, by the four point condition we have
$d(x,gy)+d(y,gx)\leq d(x,gx)+d(y,gy)+2\delta$. Therefore,
\[\begin{array}{lcl}
d(x,gx) & \geq & d(x,gy)+d(y,gx)-d(y,gy)-2\delta \\
   & = & d(x,y)+d(y,gy)-2\alpha + d(x,y)+d(y,gy)-2\beta
   -d(y,gy)-2\delta\\
   & = & 2d(x,y)+d(y,gy)-2(\alpha+\beta+\delta)\\
   & \geq & 2d(x,y)+d(y,gy)-2(3\delta+2M).
\end{array}\]

In the second case, using the four point condition again, we have
$d(x,gy)+d(y,gx)\leq d(x,y)+d(gx,gy)+2\delta$. Thus
\[d(x,y)+d(y,gy)-2\alpha+d(x,y)+d(y,gy)-2\beta\leq
d(x,y)+d(gx,gy)+2\delta\] and hence $d(y,gy)\leq \alpha
+\beta+\delta\leq 3\delta+2M$. This completes the proof.
\end{proof}

\begin{lem} \label{lem:key}Let $G$ be a relatively hyperbolic group.
Then the inner automorphism group ${\rm Inn}(G)$ of $G$ is of
finite index in ${\rm Conj}(G)$.
\end{lem}

\begin{proof} The proof of the lemma is a generalization of the proof of
\cite[Lemma 2.2]{MS}. In this case instead of the Cayley graph of
$G$ we use the $\delta$-hyperbolic metric space $X$ (in the sense
that every geodesic triangle in $X$ is $\delta$-thin) on which $G$
acts by isometries. One essential difference between the two
proofs is the existence of parabolic isometries in the case of
relatively hyperbolic groups.

Suppose on the contrary that ${\rm Inn}(G)$ is of infinite index
in ${\rm Conj}(G)$ and fix an infinite sequence
$f_{1},f_{2},\dots,f_{n},\dots$ of conjugating automorphisms of
$G$ representing pairwise distinct cosets of ${\rm Inn}(G)$ in
${\rm Conj}(G)$. In particular, $G$ is neither finite nor
virtually infinite cyclic. Let $\lambda_{i}=\inf\limits_{x\in
X}\max\limits_{s\in S} d(x,f_{i}(s)x)$, where $S$ is a fixed
finite generating set of $G$ closed under inverses, and let
$x_{i}^{0}\in X$ such that $\max\limits_{s\in
S}d(x_{i}^{0},f_{i}(s)x_{i}^{0})\leq \lambda_{i}+\frac{1}{i}$. As
shown in the proof of Theorem 1.2 in \cite{BS}, the sequence
$\lambda_{i}$ converges to infinity. Hence, for a given
non-principal ultrafilter $\omega$ on $\mathbb{N}$ the based
ultralimit $(X_{\omega},d_{\omega},x_{\omega}^0)$ of the sequence
of based metric spaces $(X,d_{i},x_{i}^{0})$, where
$d_{i}=\frac{d}{\lambda_{i}}$, is an $\mathbb{R}$-tree. Moreover,
there is an induced non-trivial isometric $G$-action on
$(X_{\omega},d_{\omega},x_{\omega}^0)$ (i.e. $G$ has no a global
fixed point in $X_{\omega}$), given by $g\cdot
(x_{i})=\big(f_{i}(g)x_{i}\big)$.

We shall show again that this action has a global fixed point.
Suppose $g$ is an element of $G$ acting as a hyperbolic isometry
on $X_{\omega}$ with translation length $\tau_{\omega}(g)$, and
fix an element $x=(x_{i})\in X_{\omega}$ on the axis of $g$. Then
\\
$\lim_{\omega}d_{i}\big(f_{i}(g)^{2}x_{i},x_{i}\big)=d_{\omega}(g^{2}x,x)=
2\tau_{\omega}(g)=2d_{\omega}(gx,x)=2\lim_{\omega}d_{i}\big(f_{i}(g)x_{i},x_{i}\big)$
and so
$\textrm{lim}_{\omega}\Big(2d_{i}\big(f_{i}(g)x_{i},x_{i}\big)-d_{i}
\big(f_{i}(g)^{2}x_{i},x_{i}\big)\Big) =0 $.

For each index $i$, let $Y_{i}=\big\{y\in X: \tau(f_{i}(g))\leq
d(y,f_{i}(g)y)\leq \tau(f_{i}(g))+\frac{1}{i}\big\}$, where
$\tau(f_{i}(g))$ is the minimal displacement of $f_{i}(g)$, and
choose $y_{i}\in Y_{i}$ such that $d(x_{i},y_{i})\leq
d(x_{i},Y_{i})+1$.

By Lemma \ref{lem:inequalities}, there are non-negative constants
$C(\delta)$ and $K(\delta)$, depending only on $\delta$, such that
$d\big(f_{i}(g)x_{i},x_{i}\big)\geq
2d(x_{i},y_{i})+d\big(f_{i}(g)y_{i},y_{i}\big)-K(\delta)$ whenever
$d(y_{i},f_{i}(g)y_{i})>C(\delta)$. Let $I$ denote the subset of
$\mathbb{N}$ consisting of those indices $i$ for which
$d(y_{i},f_{i}(g)y_{i})>C(\delta)$. The maximality of $\omega$
implies that it contains exactly one of $I$, $\mathbb{N}- I$.

We consider the two cases separately.

\textbf{Case 1:} $I\in \omega$. In this case for each $i\in I$ we
have $d\big(f_{i}(g)x_{i},x_{i}\big)\geq
2d(x_{i},y_{i})+d\big(f_{i}(g)y_{i},y_{i}\big)-K(\delta)\geq
2d(x_{i},y_{i})+\tau(f_{i}(g))-K(\delta)$. On the other hand,
$d\big(f_{i}(g)x_{i},x_{i}\big)\leq
2d(x_{i},y_{i})+d\big(f_{i}(g)y_{i},y_{i}\big)\leq
2d(x_{i},y_{i})+\tau(f_{i}(g))+\frac{1}{i}$. It follows that
$|d\big(f_{i}(g)x_{i},x_{i}\big)-\tau(f_{i}(g))|\leq
2d(x_{i},y_{i})+K(\delta)+\frac{1}{i}$. Now, it is easy to verify
that
\[2d\big(f_{i}(g)x_{i},x_{i}\big)-d\big(f_{i}(g)^{2}x_{i},x_{i}\big)\geq\]
\[4d(x_{i},y_{i})+2\tau(f_{i}(g))-2K(\delta)-[2d(x_{i},y_{i})+2d\big(f_{i}(g)y_{i},y_{i}\big)]\geq\]
\[4d(x_{i},y_{i})+2\tau(f_{i}(g))-2K(\delta)-[2d(x_{i},y_{i})+2\tau(f_{i}(g))+\frac{2}{i}]=2d(x_{i},y_{i})-2K(\delta)-\frac{2}{i},\]
and hence $2d(x_{i},y_{i})\leq
2d\big(f_{i}(g)x_{i},x_{i}\big)-d\big(f_{i}(g)^{2}x_{i},x_{i}\big)+2K(\delta)+\frac{2}{i}$.
Finally, \[\left|\tau_{\omega}(g) -
\frac{\tau(f_{i}(g))}{\lambda_{i}}\right| \leq
\left|\tau_{\omega}(g)-
d_{i}\big(f_{i}(g)x_{i},x_{i}\big)\right|+\left|
d_{i}\big(f_{i}(g)x_{i},x_{i}\big)-
\frac{\tau(f_{i}(g))}{\lambda_{i}}\right|\] \[ \leq
\left|\tau_{\omega}(g)-
d_{i}\big(f_{i}(g)x_{i},x_{i}\big)\right|+2\frac{d(x_{i},y_{i})}{\lambda_{i}}+\frac{K(\delta)}{\lambda_{i}}+\frac{1}{i\lambda_{i}}\]
\[  \leq  \left| \tau_{\omega}(g)- d_{i}(f_{i}(g)x_{i},x_{i})\right|+
2d_{i}\big(f_{i}(g)x_{i},x_{i}\big)-d_{i}\big(f_{i}(g)^{2}x_{i},x_{i}\big)+3\frac{K(\delta)}{\lambda_{i}}+\frac{3}{i\lambda_{i}},\]
for all $i\in I$.

Since the $\omega$-limit of the right-hand side of the above
inequality is $0$ and $\tau(f_{i}(g))=\tau(g)$ for all $i$,
$f_{i}(g)$ being a conjugate of $g$ for each $i$, it follows that
$\tau_{\omega}(g)=0$, which contradicts the assumption that
$\tau_{\omega}(g)$ is strictly positive.

\textbf{Case 2:} $\mathbb{N}-I\in \omega$. If
$\lim_{\omega}d_{i}(x_{i},y_{i})<\infty$, the sequence $y=(y_{i})$
is a point of $X_{\omega}$ fixed by $g$, since $0\leq
d_{\omega}(y,gy)=\lim_{\omega}d_{i}\big(y_{i},f_{i}(g)y_{i}\big)\leq
\lim_{\omega}\frac{C(\delta)}{\lambda_{i}}=0$, which contradicts
the choice of $g$. Hence $\lim_{\omega}d_{i}(x_{i},y_{i})=\infty$.
For each $i$, let $\gamma_{i}:[0,d_{i}(x_{i},y_{i})]\rightarrow
X_{i}$ be a geodesic from $x_{i}$ to $y_{i}$. Since for every
$t\geq 0$ the set of indices $i$ for which $t$ lies in the domain
of $\gamma_{i}$, is contained in $\omega$, we can define a
geodesic ray $\gamma:[0,\infty)\rightarrow X_{\omega}$ by
$\gamma(t)=\big(\gamma_{i}(t)\big)_{i}$, which is asymptotic to an
ideal point $\tilde{y}\in \partial X_{\omega}$. We will show first
that $g$ fixes $\tilde{y}$, and then that $\tilde{y}$ is not one
of the points at infinity determined by the axis of $g$,
contradicting the fact that any hyperbolic isometry of a tree
fixes exactly two points at infinity.

\begin{claim} $g$ fixes $\tilde{y}$.\end{claim}

\noindent\emph{Proof.} It suffices to show that the geodesics
$g\gamma$ and $\gamma$ are asymptotic, i.e. that
$\textrm{sup}_{t}\,d_{\omega}(g\gamma(t),\gamma(t))< \infty$. Let
$t\geq 0$. For each $i$ big enough, we consider the quadrilateral
defined by the geodesics $f_{i}(g)\gamma_{i}$, $\gamma_{i}$,
$[x_{i},f_{i}(g)(x_{i})]$ and $[y_{i},f_{i}(g)(y_{i})]$. Since the
space $X_{i}$ is hyperbolic, there is a non-negative constant
$M(\delta)$, depending only on $\delta$, such that the side
$f_{i}(g)\gamma_{i}$ is contained in the
$\frac{M(\delta)}{\lambda_{i}}$-neighborhood of the union of the
other sides. Hence, the side $f_{i}(g)\gamma_{i}$ is contained in
the $R$-neighborhood of $\gamma_{i}$, where
$R=\frac{M(\delta)}{\lambda_{i}}+d_{i}\big(x_{i},f_{i}(g)(x_{i})\big)+d_{i}\big(y_{i},f_{i}(g)(y_{i})\big)$.
Let $t'\geq 0$ be such that $\gamma_{i}(t')$ is the projection of
$f_{i}(g)\gamma_{i}(t)$ on $\gamma_{i}$. Then
$t'=d_{i}\big(\gamma_{i}(t'),\gamma_{i}(0)\big)\leq
d_{i}\big(\gamma_{i}(t'),f_{i}(g)\gamma_{i}(t)\big)+
d_{i}\big(f_{i}(g)\gamma_{i}(t),f_{i}(g)\gamma_{i}(0)\big)+d_{i}\big(f_{i}(g)\gamma_{i}(0),\gamma_{i}(0)\big)$
and thus $t'-t\leq R+d_{i}\big(f_{i}(g)x_{i},x_{i}\big)$. In the
same way, we obtain that $t-t'\leq
R+d_{i}\big(f_{i}(g)x_{i},x_{i}\big)$ and therefore $|t'-t|\leq
R+d_{i}\big(f_{i}(g)x_{i},x_{i}\big)$. Now
\[\begin{array}{lcl}
d_{i}\big(f_{i}(g)\gamma_{i}(t),\gamma_{i}(t)\big) & \leq &
d_{i}\big(f_{i}(g)\gamma_{i}(t),\gamma_{i}(t')\big)+d_{i}\big(\gamma_{i}(t'),\gamma_{i}(t)\big)\\
 & \leq & R+|t-t'|\\
 & \leq &
 2\frac{M(\delta)}{\lambda_{i}}+3d_{i}\big(x_{i},f_{i}(g)(x_{i})\big)+2d_{i}\big(y_{i},f_{i}(g)(y_{i})\big).
\end{array}\]
Taking limits, we get
$d_{\omega}\big(g\gamma(t),\gamma(t)\big)\leq 3 \tau_{\omega}(g)$.
This proves the claim.

\begin{claim} $\tilde{y}$ is not one
of the points at infinity determined by the axis of $g$.
\end{claim}

\noindent\emph{Proof.} Suppose that $\tilde{y}$ is one of the ends
of the axis $A_{g}$ of $g$. Since $X_{\omega}$ is a tree, there is
$t_{0}\geq 0$ such that $\gamma(t)\in A_{g}$ for all $t\geq
t_{0}$. The assumption that $g$ acts on $A_{g}$ as a translation
of amplitude $\tau_{\omega}(g)$ implies that either
$g\gamma(t)=\gamma\big(t+\tau_{\omega}(g)\big)$ or
$g^{-1}\gamma(t)=\gamma\big(t+\tau_{\omega}(g)\big)$ for all
$t\geq t_{0}$. Suppose that
$g\gamma(t)=\gamma\big(t+\tau_{\omega}(g)\big)$ for all $t\geq
t_{0}$ (the other case is handled similarly). Then
$\lim_{\omega}d_{i}\big(f_{i}(g)\gamma_{i}(t),\gamma_{i}(t+\tau_{\omega}(g))\big)=0$.
Fixing $t\geq t_{0}$, the geodesic $\gamma_{i}$ contains the point
$\gamma_{i}(t)+\tau_{\omega}(g)$ for each $i$ sufficiently large.
Thus
$\tau_{\omega}(g)+d_{i}\big(\gamma_{i}(t+\tau_{\omega}(g)),y_{i}\big)=d_{i}\big(\gamma_{i}(t),y_{i}\big)=
d_{i}\big(f_{i}(g)\gamma_{i}(t),f_{i}(g)y_{i}\big)\leq
d_{i}\big(f_{i}(g)\gamma_{i}(t),y_{i}\big)+d_{i}\big(y_{i},f_{i}(g)y_{i}\big)$,
and so

\[\begin{array}{lcl}
\tau_{\omega}(g) & \leq &
d_{i}\big(f_{i}(g)\gamma_{i}(t),y_{i}\big)-d_{i}\big(\gamma_{i}(t+\tau_{\omega}(g)),y_{i}\big)+d_{i}\big(y_{i},f_{i}(g)y_{i}\big)\\
 & \leq & d_{i}\big(f_{i}(g)\gamma_{i}(t),\gamma_{i}(t+\tau_{\omega}(g))\big)+d_{i}\big(y_{i},f_{i}(g)y_{i}\big).\end{array}\]
It follows that $0<\tau_{\omega}(g)\leq
\lim_{\omega}d_{i}\big(f_{i}(g)\gamma_{i}(t),\gamma_{i}(t+\tau_{\omega}(g))\big)+\lim_{\omega}d_{i}\big(y_{i},f_{i}(g)y_{i}\big)=0$.
This is a contradiction, proving the claim.

So in all cases, every element of the finitely generated group $G$
fixes some point of $X_{\omega}$. This implies that the action of
$G$ on $X_{\omega}$ has a global fixed point (see
\cite[Proposition II. 2. 15]{MoSh}), which is the desired
contradiction.
\end{proof}

\begin{rem}
It follows from the above proof that
$\tau_{\omega}(g)=\lim_{\omega}\frac{\tau(f_{i}(g))}{\lambda_{i}}$
for each $g\in G$ and each sequence $(f_{i})$ of automorphisms
representing pairwise distinct elements in ${\rm Out}(G)$. The
proof of the lemma can be simplified if one at the beginning makes
use of the hypothesis that each $f_{i}$ is a conjugating
automorphism.
\end{rem}

\begin{rem} In the proof of \cite[Lemma 2.2]{MS} the points $y_{i}$
were chosen so that $f_{i}(g)$ realizes its minimal displacement
at $y_{i}$. However, this random choice could give $y_{i}$ for
which the inequality $d(f_{i}(g)x_{i},x_{i})\geq
2d(x_{i},y_{i})+d(f_{i}(g)y_{i},y_{i})-K(\delta)$ is false. The
correct way to proceed with the proof is to choose $y_{i}$ as
above. The fact that in the case of hyperbolic groups the action
(on the Cayley graph) is free and cocompact, can be used to avoid
Case 2. Indeed, in such an action the translation lengths are
bounded away from zero and therefore for each positive number $K$
there is a positive integer $m$ such that $\tau(g^{m})>K$ for all
group elements $g$ of infinite order. Thus, by replacing each
element by its $m$-th power, we can suppose that the translation
length of any $f_{i}(g)$ is big enough.
\end{rem}

\begin{proof}[Proof of Theorem \ref{thm:1}]
We consider the short exact sequence (\ref{eq:1}). Since the first
term is finite, each torsion-free subgroup of ${\rm Out}(G)$
embeds in ${\rm Aut}(G)/{\rm Conj}(G)$, which is residually finite
by \cite[Lemma 2.1]{MS}. Hence, each torsion-free subgroup of
${\rm Out}(G)$ is residually finite, from which the theorem
follows.
\end{proof}

\begin{rem} Let $\textrm{Aut}_{n}(G)$ be the subgroup of
$\textrm{Aut}(G)$ consisting of automorphisms which fix setwise
every normal subgroup of $G$. Obviously, ${\rm Conj}(G)\subseteq
\textrm{Aut}_{n}(G)$. After this paper appeared as a preprint,
Minasyan and Osin \cite{MO} proved (using completely different
methods than ours) that for any relatively hyperbolic group $G$,
${\rm Inn}(G)$ has finite index in $\textrm{Aut}_{n}(G)$ and hence
in ${\rm Conj}(G)$. Moreover, if $G$ is non-elementary and has no
non-trivial finite normal subgroups, then
$\textrm{Aut}_{n}(G)={\rm Inn}(G)$. Thus, in this case, the
hypothesis of virtual torsion freeness can be removed in Theorem
\ref{thm:1}.\end{rem}

To prove Theorem \ref{thm:2}, we need the following result.

\begin{thm}[\textrm{\cite[Corollary 5.3]{GL}}]\label{thm:GL}
Let $G=G_{1}\ast \cdots \ast G_{n}\ast F_{k}$, where each $G_{i}$
is finitely generated, freely indecomposable and not infinite
cyclic, and $F_{k}$ is a free group of rank $k$, with $n+k\geq 2$.
Suppose that each factor $G_{i}$ contains a torsion-free, normal
subgroup of finite index $H_{i}$ such that ${\rm Out}(H_{i})$ is
virtually torsion-free and the quotient $H_{i}/Z(H_{i})$ of
$H_{i}$ by its center $Z(H_{i})$ is torsion-free. Then ${\rm
Out}(G)$ is virtually torsion-free. \qed \end{thm}

We also need the following simple lemma, whose proof is left to
the reader.

\begin{lem} \label{lem:1} Every polycyclic-by-finite group $G$ has a
normal, torsion-free, finite index subgroup $H$ such that the
quotient $H/Z(H)$ of $H$ by its center $Z(H)$ is also
torsion-free. \qed \end{lem}

\begin{proof}[Proof of Theorem \ref{thm:2}]
By Dyer's results \cite{Dy}, the class of conjugacy separable
groups is closed under finite graphs of groups with finite edge
groups. Since polycyclic-by-finite groups are conjugacy separable
\cite{Fo}, it follows that $G$ is conjugacy separable.

It remains to show that ${\rm Out}(G)$ is virtually torsion-free.
By \cite[Lemma 2.4]{MS}, it suffices to find a normal subgroup $N$
of finite index in $G$ with trivial center and virtually
torsion-free outer automorphism group. Since $G$ is conjugacy
separable (and so residually finite), it has a normal subgroup of
finite index $N$ which intersects each edge group trivially. This
means that $N$ acts non-trivially on the corresponding tree of $G$
with trivial edge stabilizers and therefore $N$ admits a
non-trivial free product decomposition $N_{1}\ast \cdots \ast
N_{k}$ into freely indecomposable, polycyclic-by-finite factors.
In particular, the center of $N$ is trivial. To see that ${\rm
Out}(N)$ is virtually torsion-free, note first that the outer
automorphism group of a polycyclic-by-finite group is virtually
torsion-free being finitely generated and isomorphic to a subgroup
of $GL_{n}(\mathbb{Z})$, for some positive integer $n$ (see
\cite{We}). Thus, in view of Lemma \ref{lem:1}, the hypotheses of
Theorem \ref{thm:GL} are satisfied and so ${\rm Out}(N)$ is
virtually torsion-free.
\end{proof}

\noindent{\footnotesize \textbf{Acknowledgement.} We thank the
referee of the previous version for pointing out a gap in the
proof of Lemma \ref{lem:key}}.

\noindent Department of Mathematics, University of the Aegean, Karlovassi, 832 00 Samos,
Greece {\it email: vmet@aegean.gr}

\bigskip

\noindent Department  of Mathematics and Statistics, University of
Cyprus, P.O. Box 20537, 1678 Nicosia, Cyprus {\it email:
msikiot@ucy.ac.cy, msykiot@math.uoa.gr}
\end{document}